\documentclass[11pt,oneside]{article}
\usepackage[utf8]{inputenc}

\usepackage[english]{babel}%
\usepackage[a4paper]{geometry}
\usepackage{indentfirst}

\usepackage{amsmath,amssymb,amsfonts,amsthm,mathrsfs,stmaryrd}  
\usepackage{graphicx,caption,float}
\usepackage[usenames, dvipsnames]{xcolor}

\usepackage[hang,flushmargin]{footmisc} %%flushmargin for footnote

\usepackage[colorlinks=true, pdfstartview=FitV,linkcolor=ForestGreen,citecolor=ForestGreen, urlcolor=blue]{hyperref}

\author{
Nicolae Mihalache$^{1}$
\quad and \quad
Fran\c{c}ois Vigneron$^{2}$
}

\newcommand{\authornicu}{{\small
\textbf{Nicolae Mihalache.}
Université Paris-Est Creteil, UMR 8050 CNRS, LAMA, F-94010 Creteil, France \&
Université Gustave Eiffel, LAMA, F-77447 Marne-la-Vallée, France\\
\phantom{x}\hfill\texttt{nicolae.mihalache@u-pec.fr}
}}

\newcommand{\authorfv}{{\small
\textbf{François Vigneron.}
Universit\'{e} de Reims Champagne-Ardenne,
Laboratoire de Mathématiques de Reims,
UMR 9008 CNRS, Moulin de la Housse, BP 1039,
F-51687 Reims\\
\phantom{x}\hfill\texttt{francois.vigneron@univ-reims.fr}
}}

\newtheorem{thm}{Theorem}
\newtheorem{prop}{Proposition}

\newtheorem{remark}[prop]{Remark}
\newtheorem{example}[prop]{Example}

\newcommand{\N}{\mathbb{N}}
\newcommand{\Z}{\mathbb{Z}}

\newcommand{\C}{\mathbb{C}}
\newcommand{\CC}{\overline{\mathbb C}}

\newcommand\ii[2]{\ensuremath{\llbracket{#1}, {#2} \rrbracket}}
\renewcommand{\div}{\operatorname{Div}}

\newcommand{\M}{\mathcal{M}}
\newcommand{\julia}{\mathcal{J}}
\newcommand{\fatou}{\mathcal{F}}
\DeclareMathOperator{\hyp}{Hyp}
\DeclareMathOperator{\mis}{Mis}
\let\oldRe=\Re
\renewcommand{\Re}{\operatorname{\oldRe e}}

\newcommand{\eg}[1][~]{\textit{e.g.}#1}
\newcommand{\ie}[1][~]{\textit{i.e.}#1}

\newcommand{\defequal}{=}

\begin{document}

\title{Factorization of the quadratic\\
Misiurewicz-Thurston polynomials}
\maketitle

\begin{abstract}  
This note provides the complete factorization of the Misiurewicz-Thurston polynomial $q_{\ell,n}$ over $\C$,
which plays a central role in the study of the Mandelbrot set.
The roots can be classified into two categories.
First, there are hyperbolic points $\hyp(k)$ for any divisor $k$ of $n$, which are parameters whose critical orbits are of exact period~$k$. Those are roots of $q_{\ell,n}$ with multiplicity $\left\lfloor \frac{\ell -1}{k} \right\rfloor + 2$.
Next are the points $\mis(j,k)$ for $2\leq j\leq \ell$ whose critical orbits are pre-periodic of
exact period $k$ with an exact pre-period $j$. Those are simple roots of $q_{\ell,n}$.
\end{abstract}

%%%%%%%%%%%%%%%%%%%%%%%%%%%%%%%%

In this article, we are interested in the factorization of two important families of polynomials.
The first family is defined recursively by
\begin{equation}\label{def:hyp_poly}
p_0(z) \defequal 0, \qquad p_{n+1}(z) \defequal p_n(z)^2 + z.
\end{equation}
For $n\geq 1$, one has $\deg p_n = 2^{n-1}$.
The second family is a two parameters family defined for $\ell,n\in\N$ by %of \textit{Misiurewicz-Thurston polynomials} 
\begin{equation}\label{def:mis_poly}
q_{\ell,n}(z)=p_{\ell+n}(z) - p_\ell(z).
\end{equation}
The factorisation of $p_n$ over $\C$ is well known \cite{DOUHUB82} and is recalled first to fix essential notations.
The main result of this article is the proof of Theorem~\ref{thm:fact_mis} below, which gives the complete
factorization of $q_{\ell,n}$ in terms of simple factors whose dynamical significance is explained
briefly in the first section.
This result plays a role in the systematic computation of hyperbolic centers
and of pre-periodic parameters in the Mandelbrot set, up to period $n=41$ for
$p_n$ and all types $\ell+n\leq 35$ for $q_{\ell,n}$, as detailed in~\cite{MV2025}, yet remains of general interest.

\section{About the Mandelbrot set}
\label{par:mandel}

\begin{figure}[H]
%\captionsetup{width=.9\linewidth}
\begin{center}
\includegraphics[width=.55\textwidth]{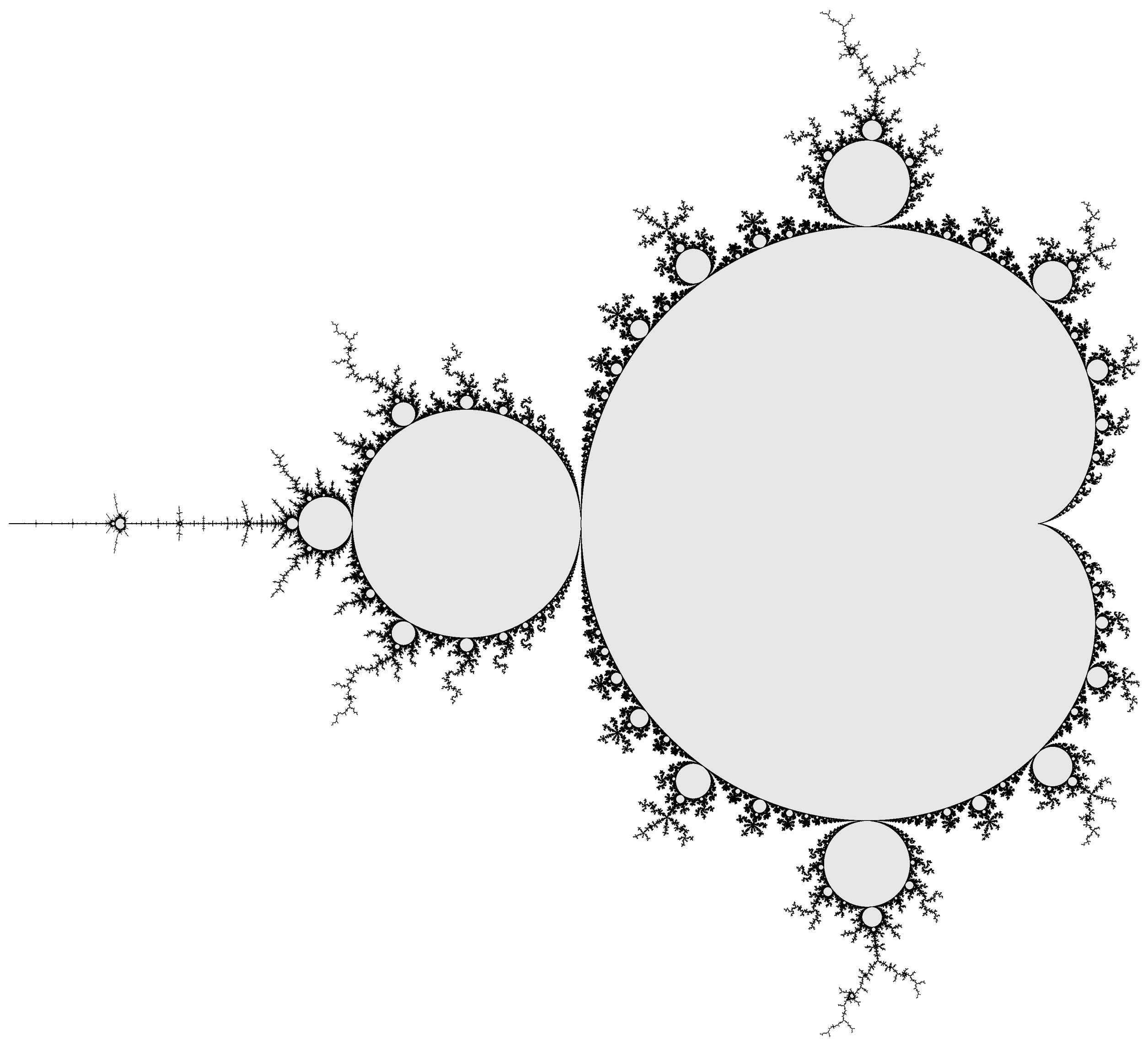}
\caption{\label{fig:mandelbrot}\small\sl
The Mandelbrot set~$\M$ with $\partial\M$ in black and in gray, the interior of $\M$.}% caption
\end{center}
\end{figure}

The polynomials $p_n$ and $q_{\ell,n}$ are closely related to the \textit{Mandelbrot} set $\M$ (Fig.~\ref{fig:mandelbrot}).
This set is composed of the parameters~$c\in\C$ for which the sequence~$(p_n(c))_{n\in\N}$
remains bounded.
This sequence is the orbit by iterated compositions of the only critical point $z=0$ of the map $f_c(z) \defequal z^2 +c$, \ie
\begin{equation}\label{eq:defFundMap}
\forall n\in\N,\qquad
f_c^n (0)\defequal \underset{\text{$n$ times}}{\underbrace{\,f_c \circ \ldots\circ f_c}}(0)  = p_n(c).
\end{equation}

For $c\in\C$, the dynamics of $f_c$ splits $\CC=\C\cup\{\infty\}$ in two complementary sets.
The \textit{Fatou} set~$\fatou_c$ is the open subset composed by the points $z\in\CC$
in the neighborhood of which, the sequence $(f_c^n)_{n\in\N}$ is a normal family
\ie precompact in the topology of local uniform convergence.
On the contrary, on the \textit{Julia} set~$\julia_c=\CC\backslash\fatou_c$, the dynamics is chaotic.
Both~$\fatou_c$ and~$\julia_c$ are fully invariant (\ie invariant sets of the forward and backwards dynamics)
and that $c\in\M$ if and only if $\julia_c$ is connected.
For a review of the properties of Fatou and Julia sets, see~\eg\cite{CARLESON}, \cite{MILNOR} for polynomial and rational maps
and~\cite{Ber93}, \cite{MRW22} for entire and meromorphic functions.

\subsection{Hyperbolic centers}

For $n\geq1$, the roots of $p_n$ are parameters $c\in\M$, called \textit{hyperbolic centers},
whose critical orbits are periodic of period $n$.
The roots of $p_n$ are simple~\cite{DOUHUB82, BUFF2018}
and the polynomial $p_n(z)$ is divisible by $p_k(z)$ for any divisor $k$ of $n$.
A.~Douady and J.H.~Hubbard \cite{DOUHUB82, DOUHUB84} have shown that the hyperbolic centers are interior points of~$\M$,
with at most one hyperbolic center per connected component of the interior.
The set of all hyperbolic centers is also dense in the boundary of $\M$ in the sense that
its closure in $\C$ contains $\partial\M$.

\medskip
Let us define the set of hyperbolic centers of \textit{order} $n\geq1$
as the subset of $p_n^{-1}(0)$ whose minimal (or fundamental) period is exactly $n$; in other terms:
\begin{equation}\label{def:hyp}
\hyp(n) \defequal \left\{ z\in p_n^{-1}(0) \,\big\vert\, \forall k\in\div(n)^\ast, \: p_k(z)\neq 0 \right\}
\end{equation}
where $\div(n) \defequal \{ k\in\N^\ast \,;\, \exists k'\in\N, \enspace \enspace k k' = n\}$
and $\div(n)^\ast=\div(n)\backslash\{n\}$ is the set of strict divisors of $n$.
For example,~$\hyp(1) = \{0\}$ and $\hyp(2) = \{-1\}$. The reduced polynomial, also known as Gleason's polynomial, is:
\begin{equation}\label{eq:hyp_red}
h_n(z) \defequal \prod\limits_{r\in\hyp(n)}(z-r).
\end{equation}
The polynomials $p_n$ and $h_n$ have integer coefficients.
While expected, the irreducibility of $h_n$ over $\Z[z]$ remains conjectural; see \cite[last remark of~\S3]{HT15} and \cite[p.155]{SS2017}.

\begin{thm}[\cite{DOUHUB82}]\label{thm:countingHyp}
The complete factorization of $p_n$ is
\begin{equation}\label{eq:HypFactor}
p_n(z) = \prod_{k\vert n}h_k(z).
\end{equation}
Moreover, the cardinal of $\hyp(n)$ is given by
\begin{equation}\label{eq:HypCount}
\left|\hyp(n)\right| = \sum_{k\vert n} \mu(n/k)2^{k-1}
\end{equation} 
where $\mu$ is the M\"obius function, \ie
\[
\mu(n)=\begin{cases}
(-1)^{\nu} & \text{if $n$ is square free and has $\nu$ distinct prime factors,}\\
0 & \text{if $n$ is not square free.}
\end{cases}
\]
\end{thm}
\noindent
The notation $k\vert n$ means $k\in\div(n)$ and $|S|$ denotes
the cardinal of a finite set.
The identity~\eqref{eq:HypFactor} is well known and is key in the count of hyperbolic centers.
The Online Encyclopedia of Integer Sequences~\cite{CARDHk} attributes~\eqref{eq:HypCount} to 
Warren D. Smith and Robert Munafo (2000), sadly with no published reference.
For a more general result that applies to iterated maps of the form $z^d+c$, see~\cite{HT15}.

\begin{figure}[H]
\captionsetup{width=.95\linewidth}
\begin{center}
\includegraphics[width=\textwidth]{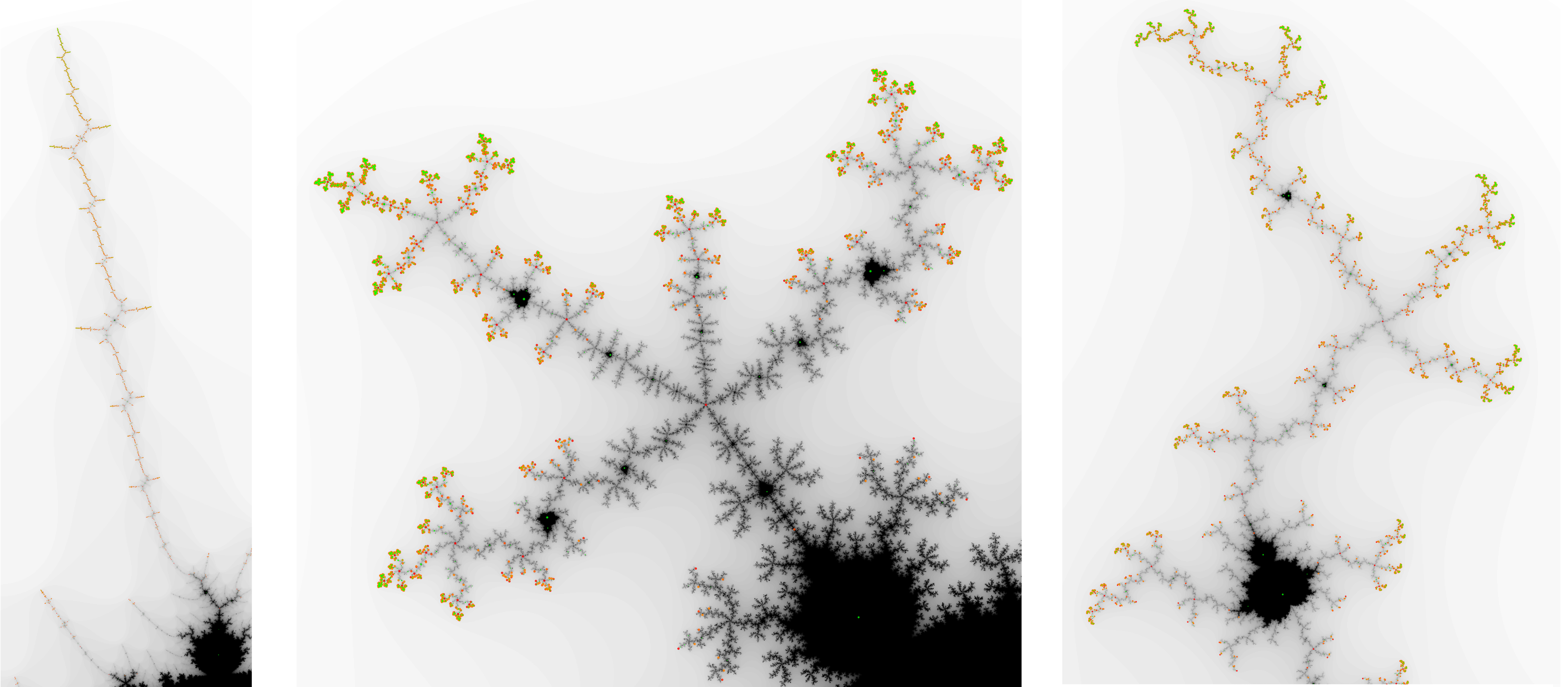}
\caption{\label{fig:HypMis}
The hyperbolic points $\hyp(n)$ for $n\leq18$ in green and
the Misiurewicz-Thurston parameters $M_{\ell,n}=\mis(\ell,n)$ with $\ell+n\leq16$ in red
in different parts of the Mandelbrot set~$\M$, with $\Re c$ increasing from left to right
between images.
}% caption
\end{center}
\end{figure}

\subsection{Pre-periodic or Misiurewicz-Thurston parameters}

For integers $n\geq1$ and $\ell\geq2$, the roots of $q_{\ell,n}$ (defined in \eqref{def:mis_poly}) are parameters $c\in\partial\M$, called \textit{pre-periodic parameters}
or \textit{Misiurewicz-Thurston} points, whose critical orbits are pre-periodic, that is it becomes periodic of period $n$ after the first $\ell$ steps.
For $\ell\in\{0,1\}$, one can check immediately that
\begin{equation}\label{eq:hypCase}
q_{0,n}(z) = p_n(z) \qquad\text{and}\qquad q_{1,n}(z) = p_n^2(z).
\end{equation}
The dynamical reason for the second identity in~\eqref{eq:hypCase} is that $0$ is the only pre-image of~$c$ under $f_c$
so a pre-periodic orbit of type $(1,n)$ starting at zero actually loops back to zero, so is periodic of period $n$.
The polynomial $q_{\ell,n}(z)$ is divisible by $q_{\ell,k}(z)$ for any divisor $k$ of $n$ and by $q_{\ell',n}(z)$ for any $\ell'\leq\ell$.
Some of the roots of $q_{\ell,n}(z)$ have multiplicity; however, the polynomial
\begin{equation}\label{eq:polySimple}
s_{\ell,n}(z)=\frac{q_{\ell,n}(z)}{q_{\ell-1,n}(z)} \in \Z[z]
\end{equation}
has simple roots \cite{BUFF2018}, \cite[Lemma~3.1]{HT15}. 
Douady-Hubbard \cite{DOUHUB82, DOUHUB84} have shown that the set of all pre-periodic points is dense in the boundary of $\M$.
Visually, those points are either branch tips, centers of spirals or points where branches meet (see Fig.~\ref{fig:HypMis}).

\medskip
By analogy with the hyperbolic case, let us define the set of \textit{Misiurewicz points of type $(\ell,n)$}
as the subset of $q_{\ell,n}^{-1}(0)$ whose dynamical parameters are exactly $\ell$ and $n$; in other terms:
\begin{equation}\label{def:mis}
\mis(\ell,n) \defequal \left\{ z\in q_{\ell,n}^{-1}(0)
\left|\begin{array}{l}
q_{\ell-1,n}(z)\neq 0,\\[3pt]
\forall k\in\div(n)^\ast, \: q_{\ell,k}(z)\neq 0
\end{array}\right.\right\}.
\end{equation}
We call $\ell+n$ the \textit{order} of $\mis(\ell,n)$ because $q_{\ell,n}$ is a polynomial of degree $2^{\ell+n-1}$.
The reduced polynomial whose roots are exactly $\mis(\ell,n)$ is denoted by
\begin{equation}\label{eq:mis_red}
m_{\ell,n}(z) \defequal \prod\limits_{r\in\mis(\ell,n)}(z-r) \in \Z[z].
\end{equation}
Note that $\mis(0,n)=\hyp(n)$ and $\mis(1,n)=\emptyset$ because of~\eqref{eq:hypCase}.
The count of the Misiurewicz points has been established by B.~Hutz and A.~Towsley \cite[Cor.~3.3]{HT15}:
\begin{equation}\label{eq:MisCount}
\left|\mis(\ell,n)\right| = \Phi(\ell,n)\left|\hyp(n)\right|
\end{equation}
where
\[
\Phi(\ell,n) = \begin{cases}
1 & \text{if } \ell = 0,\\
2^{\ell-1} -1 & \text{if } \ell \neq 0 \text{ and }n\vert \ell-1,\\
2^{\ell-1} & \text{otherwise.}
\end{cases}
\]

\medskip
To get a complete factorization of $q_{\ell,n}(z)$ in terms of~\eqref{eq:hyp_red} and~\eqref{eq:mis_red},
one needs to understand the multiplicity of the hyperbolic factors. We claim the following result.

\begin{thm}\label{thm:fact_mis}
For $\ell\in\N$ and $n\in \N^\ast$, one has
\begin{equation}\label{eq:MisFactor}
q_{\ell,n}(z) = \prod_{k\vert n} \left( h_k(z)^{\eta_\ell(k)} \prod_{j=2}^{\ell} m_{j,k}(z) \right)
\end{equation}
where the multiplicity $\eta_\ell(k)$ is given by\footnote{In equation~\eqref{eq:defMultHyp}, the notation $\lfloor\cdot\rfloor$ represents the floor function, \ie $\lfloor x \rfloor = n$ if $n\in\Z$ and $x\in [n,n+1)$.}
\begin{equation}\label{eq:defMultHyp}
\eta_\ell(k) \defequal \left\lfloor \frac{\ell -1}{k} \right\rfloor + 2.
\end{equation}
In other words, the roots of $q_{\ell,n}$ are composed of the points $\hyp(k)$ for any divisor $k$ of $n$,
which are roots with multiplicity $\eta_\ell(k)$, and of the points $\mis(j,k)$ for $2\leq j\leq \ell$, which are simple roots.
\end{thm}
\noindent
This result appears to be new. We propose a direct proof in Section~\ref{par:factorization}.
In \cite{MV2025}, we use this result to compare the performance of a new splitting algorithm in the case of simple roots, like
for $p_n$, or in the presence of high-multiplicity roots, like for $q_{\ell,n}$. Knowing the exact multiplicity is crucial for the
interpretation of the numerical benchmarks.

\section{Proof of Theorem~\ref{thm:fact_mis}}\label{par:factorization}

Let us give here a direct proof of the factorization theorem.

\begin{proof}
As mentioned above, thanks to \cite{HT15}, only the multiplicity of the hyperbolic factors has to be established.
For $\ell\in\{0,1\}$, the formula~\eqref{eq:MisFactor} boils down to~\eqref{eq:hypCase} so
we will suppose $\ell\geq 2$ from now on. In particular, $\eta_\ell(k)\geq2$.

\medskip
First, let us check that the multiplicities in~\eqref{eq:MisFactor} are consistent with $\deg q_{\ell,n}$, \ie:
\begin{equation}\label{eq:MisDeg}
\sum_{k\vert n} \left( \eta_\ell(k)\left|\hyp(k)\right| + \sum_{j=2}^{\ell} \left|\mis(j,k)\right| \right) = 2^{\ell+n-1}.
\end{equation}
Indeed, using~\eqref{eq:HypCount} and~\eqref{eq:MisCount} and a geometric sum, the left-hand side equals
\begin{align*}
\sum_{k\vert n} & \left|\hyp(k)\right|\left( \eta_\ell(k) + \sum_{j=2}^{\ell} \Phi(j,k) \right) \\ & =
\sum_{k\vert n} \left(\sum_{m\vert k} \mu(k/m)2^{m-1}\right)
\left( 2^\ell + \left\lfloor \frac{\ell -1}{k} \right\rfloor - \sum_{j=2}^{\ell}\delta_{k\vert j-1} \right)
\end{align*}
where $\delta_{k\vert j-1} =1$ if $k$ divides $j-1$ and $0$ otherwise.
For any integers $\lambda, k \in \N^\ast$, let us observe that 
\begin{equation}\label{eq:floorSumDelta}
\left\lfloor \frac{\lambda}{k} \right\rfloor = 
\left| \left\{ j\in\ii{1}{\lambda} \,;\,  k \vert j\right\} \right|
= \sum_{j=1}^{\lambda}\delta_{k\vert j}.
\end{equation}
The claim~\eqref{eq:MisDeg} thus follows from Möbius inversion formula \cite{Mob1832}-\cite{Rota1963}:
\begin{equation}
\sum_{k\vert n} \left(\sum_{m\vert k} \mu(k/m)2^{m}\right)= 2^{n}.
\end{equation}

\medskip
The case of $\hyp(1)=\{0\}$
is essentially based on the well known fact (proven by direct induction)
that the trailing coefficients of $p_n(z)$ stabilize as $n$ grows; more precisely:
\begin{equation}
\forall k\in\N, \qquad p_n(z) \equiv p_{n+k}(z) \mod z^{n+1},
\end{equation}
thus $q_{\ell,n}(z)\equiv 0 \mod z^{\ell+1}$. On the other hand, using~\eqref{eq:factorDiffSquares} below and
$p_k(z) = z \mod z^2$ for $k\geq1$, one has by induction:
\begin{equation}
q_{\ell,n}(z)\equiv 2^{\ell-1} z^{\ell+1} \mod z^{\ell+2},
\end{equation}
so zero is exactly of multiplicity $\ell+1 = \eta_\ell(1)$.

\medskip
Let us now consider the case of $\hyp(k)$ for $k\vert n$ and $1<k\leq n$.
Note that if $n$ is prime, there is only one hyperbolic factor left, namely $k=n$ and
the factorization~\eqref{eq:MisFactor} follows for example from an argument of
divisibility and the identity of degrees~\eqref{eq:MisDeg}.
We can however treat the general case in a unified way, regardless of wether $n$ is composite or not
and prove~\eqref{eq:MisFactor} by recurrence on $\ell\geq 2$. One has
\begin{align}
q_{\ell,n}(z)
&= p_{\ell+n}(z) - p_\ell(z)  \notag\\
&= p_{\ell+n-1}^2(z) -p_{\ell-1}^2(z)  \notag\\
&= q_{\ell-1,n}(z) (p_{\ell+n-1}(z)+p_{\ell-1}(z)) \label{eq:factorDiffSquares}\\
&= q_{\ell-1,n}(z) (q_{\ell-1,n}(z)+2p_{\ell-1}(z)). \notag
\end{align}
The recurrence assumption reads
\begin{equation}\label{eq:proofA}
q_{\ell-1,n} = \prod_{k\vert n} \left( h_k^{\eta_{\ell-1}(k)} \prod_{j=2}^{\ell-1} m_{j,k} \right)
\end{equation}
with $\eta_{\ell-1}(k)\geq 2$ and, from the hyperbolic case, we know that
\[
p_{\ell-1} = \prod_{j\vert \ell-1} h_j.
\]
If $k$ divides $n$ but not $\ell-1$, then $h_k$ divides $q_{\ell-1,n}(z)$ but is relatively prime to $p_{\ell-1}$ in $\C[z]$,
so $h_k$ and $q_{\ell-1,n}+2p_{\ell-1}$ are also relatively prime.
Alternatively, if $k$ divides both $n$ and $\ell-1$ then $h_k$ divides $q_{\ell-1,n}+2p_{\ell-1}$.
%However $h_k^2$ does not because it is a factor of $q_{\ell-1,n}$ (see the Lemma~\ref{LM} below) but  not of~$p_{\ell-1}$.
Finally, the polynomials $m_{\ell,k}$ for $k\vert n$ are known factors of $q_{\ell,n}$ that are relatively prime to all $h_j$ and
to $q_{\ell-1,n}$, and
as $q_{\ell,n}/q_{\ell-1,n}$ has simple roots, they are simple factors of $q_{\ell,n}$. There are no other
Misiurewicz-type factors and all hyperbolic factors are relatively prime with one another. We have thus established that:
\begin{equation}\label{eq:proofB}
  \left(\prod_{k\vert \operatorname{gcd}(n,\ell-1)} h_k \right)
\left(\prod_{k\vert n} m_{\ell,k} \right) q_{\ell-1,n} \quad\text{divides}\quad q_{\ell,n}.
\end{equation}

Let us observe that~\eqref{eq:defMultHyp} and \eqref{eq:floorSumDelta} imply that, for $\ell, k\geq2$:
\begin{equation}\label{eq:recEta}
\eta_{\ell-1}(k) =
\begin{cases}
\eta_\ell(k)-1 & \text{if }k\vert \ell-1,\\
\eta_\ell(k) & \text{otherwise.}
\end{cases}
\end{equation}
Combining \eqref{eq:proofA}, \eqref{eq:proofB} and~\eqref{eq:recEta}, we get that
\[
\prod_{k\vert n} \left( h_k^{\eta_{\ell}(k)} \prod_{j=2}^{\ell} m_{j,k} \right) \quad\text{divides}\quad q_{\ell,n}.
\]
Finally the identity of degrees~\eqref{eq:MisDeg} and the unitary assumption of $q_{\ell,n}$
ensure that~\eqref{eq:MisFactor} holds for the next pre-period.
\end{proof}

\begin{remark}
The previous proof also establishes the following identity regarding~\eqref{eq:polySimple}:
\begin{equation}\label{eq:polyMisSimple}
s_{\ell,n}=\frac{q_{\ell,n}}{q_{\ell-1,n}} = p_{\ell+n-1}+p_{\ell-1} = \left(\prod_{k\vert \operatorname{gcd}(n,\ell-1)} h_k \right)
\left(\prod_{k\vert n} m_{\ell,k} \right).
\end{equation}
\end{remark}

\begin{example}
To illustrate our method, the simplest case of a composite period is:
\[
q_{2,4}(z) = p_6(z) - p_2(z) = p_5^2(z) -p_1^2(z) = q_{1,4}(z) (p_5(z)+p_1(z)) = p_4^2(z) (p_4^2(z)+2z).
\]
As $p_4 = h_1 h_2 h_4$, the polynomial
$q_{2,4}$ is divisible by $h_4^2$ but not by $h_4^3$. The exponent of $\hyp(4)$ is therefore exactly 2
and we are left with a single hyperbolic factor, namely~$\hyp(2)$. Using the degree identity~\eqref{eq:MisDeg}
thus ensures that
\[
q_{2,4}(z) = h_1^3(z)h_2^2(z)h_4^2(z) m_{2,1}(z) m_{2,2}(z) m_{2,4}(z)
\]
and
\[
s_{2,4}(z) = p_4^2(z)+2z = h_1(z) m_{2,1}(z) m_{2,2}(z) m_{2,4}(z).
\]
\end{example}

%%%%%%%%%%%%%%%%%%%%%%%%%%%%%%%%
\bibliographystyle{alpha}\small
\newcommand\OFISlabel{}
\newcommand\OF[2]{OEIS}

\vspace*{3em}\noindent
$^{\small 1}$ \authornicu\\[2ex]
$^{\small 2}$ \authorfv

\end{document}